\documentstyle[12pt]{article}

\def\mcc{M\raise.5ex\hbox{c}C}
\def\mccarthy{M\raise.5ex\hbox{c}Carthy}

\def\ltm{L^2(\mu)}

\def\ptm{P^2(\mu)}

\def\hs{\htm^*}

\def\N{N^{+}}

\def\tm{T_{\bar m}}

\def\i{\infty}

\def\V{\Vert}

\def\la{Lat({\cal A})}

\def\hs{{\cal H}_s}

\def\H{\hbox{\xa H}}

\input amssym.def
\input amssym.tex
\newcommand{\dN}{{\Bbb N}}

\def\la{\langle}
\def\ra{\rangle}
\def\H{\cal H}
\def\N{\dN}
\def\fc{F_{c,w}}

\def\zs{(z_1,\dots,z_d)}
\def\cd{{\Bbb C}^d}

\def\bd{B_d}
\def\nb{N^+(\bd)}

\def\hbd{H^2({B_d})}

\def\hs{{\H}_{s}}

\def\tmh{T_{\bar m}^{H^2({\bd})}}
\def\tm{T_{\bar m}}

\def\hibd{H^\i({\bd})}
\def\hid{H^\i({B_1})}

\advance\topmargin by -.6in
\advance\textheight by 1.2in
\advance\oddsidemargin by -.5in
\advance\textwidth by 1.05in

\newtheorem{theorem}{Theorem}[section]

\newtheorem{lemma}[theorem]{Lemma}

\begin{document}
\setlength{\baselineskip}{21pt}
\title{ Moduli of bounded holomorphic functions in the ball}
\author{Boris Korenblum
\\{ State University of New York, Albany, New York 12222, U.S.A.}
\\{John E. M\raise.5ex\hbox{c}Carthy
\thanks{
Partially supported by the National Science Foundation
grant DMS 9296099}}
\\{ Washington University, St. Louis, Missouri 63130, U.S.A.}
\date{February 25 1993}
}

\bibliographystyle{plain}
\maketitle
\begin{abstract}
We prove that there is a continuous
non-negative function $g$ on the unit sphere in $\cd$, $d \geq 2$,
whose logarithm is integrable with respect to Lebesgue measure,
and which vanishes at only one point, but such that no non-zero
bounded analytic function $m$ in the unit ball, with boundary
values $m^\star$, has $|m^\star| \leq g$ almost everywhere.
The proof analyzes the common range of co-analytic Toeplitz 
operators in the Hardy space of the ball.
\end{abstract}
\baselineskip = 18pt
\setcounter{section}{-1}

\section{Introduction}
Let $B_d$ be the unit ball in $\cd$, $S_d$ be the boundary
of $B_d$, and $\sigma_d$ be normalized Lebesgue measure on $S_d$.
Every function $m$ in $H^\i(B_d)$, the space of bounded analytic
functions in $B_d$, has radial limits $\sigma_d$-almost everywhere
on $S_d$, defining a function $m^\star$ on the sphere; conversely,
$m$ can be recovered from $m^\star$ by integrating against the
Szeg\"o kernel. The problem which this paper addresses is 
when a given non-negative bounded function $g$ on $S_d$ is the
modulus $|m^\star|$ of some function $m$ in $H^\i(B_d)$.
\par
When $d=1$, the problem was completely solved by Szeg\"o
\cite{sze21}: a necessary and sufficient condition that $g$
be the modulus of a non-zero function in $\hid$ is
\setcounter{equation}{\value{theorem}}
\begin{equation}
\int_{S_d} log(g) d\sigma_d > -\i 
\label{eq:sz}
\end{equation}
\addtocounter{theorem}{1}
\par
For $d > 1$, condition~(\ref{eq:sz}) is necessary
and sufficient for $g$ to be the modulus of 
a function in the larger Nevanlinna
class $N(B_d)$, consisting of those holomorphic functions $f$
on the ball for which 
$$
T(f,1) := \sup_{0<r<1} \int_{S_d} \log^+ |f(r\zeta)| d\sigma_d(\zeta) < \i 
$$
\cite[Theorem 10.11]{rud86}. It is no longer sufficient, however, 
for $g$ to be the modulus of a bounded analytic function, 
because the function
$$
\zeta \mapsto {\rm ess\ }\sup_{-\pi \leq \theta \leq \pi}
|m^\star (e^{i\theta} \zeta) |
$$
must be lower semi-continuous on $S_d$ if $m$
is in $H^\i(B_d)$ \cite{rud86}.
In \cite[Theorem 12.5]{rud86}, 
Rudin proves that if $g$ is log-integrable,
and there exists some non-zero $f$ in $H^\i(B_d)$
with $g \geq |f^\star|$ a.e. and $g/|f^\star|$ lower semi-continuous,
then there does exist $m$ in $H^\i(B_d)$ with $g = |m^\star|$ a.e.
\par
The main result of this paper is the following:
\par
{\bf Theorem} \ref{thm:main} {\it
Let $d \geq 2$. There is a non-negative continuous function $g$ on
$S_d$, with $\int_{S_d} log (g) d \sigma_d > -\i$,
and which vanishes at only one point, but such that for no
non-zero function $m$ in $H^\i(B_d)$ is $|m^\star | \leq g$ almost
everywhere [$\sigma_d$].
}
\vskip 5pt
\par
The proof involves the analysis of co-analytic Toeplitz operators. 
If $\mu$ is a compactly supported  measure on $\cd$, 
let $\ptm$ denote the closure of the polynomials in $\ltm$,
and let $P$ denote the orthogonal projection from $\ltm$
onto $\ptm$.
If $m$ is a bounded analytic function on the support
of $\mu$, the co-analytic Toeplitz operator $T_{\bar m}^{\ptm}$ is defined
by
$$
T_{\bar m}^{\ptm} f = P \bar m f .
$$
When $\mu$ is $\sigma_d$, the space $\ptm$ is
called the Hardy space $H^2(\bd)$.
\par
The idea of the proof is as follows. 
A function $f$  is in the range of the co-analytic Toeplitz operator
$\tmh$ if and only if the linear map 
$$
\Gamma: p \mapsto \la p , f \ra_{\hbd} 
$$
is bounded on $P^2(|m^\star |^2 \sigma)$. So $f$ is in the range
of all co-analytic Toeplitz operators if and only if $\Gamma$
is a continuous linear functional on
\setcounter{equation}{\value{theorem}}
\begin{equation}
\cup_{m \in \hibd} P^2(|m^\star |^2 \sigma_d) .
\label{eq:un}
\end{equation}
\addtocounter{theorem}{1}
The Smirnov class $N^+ (\bd)$ consists of
those functions $f$ in $N(\bd)$ for which $\{ \log^+|f(r\, \cdot)| :
0<r<1 \}$ is a uniformly integrable family on $S_d$. Equipped
with the metric $\rho(f,g) = \int_{S_d} \log(1+ |f-g|) d\sigma_d$,
it becomes a topological vector space that is not locally
convex.
When $d=1$, (\ref{eq:un}) coincides (as a set) with the Smirnov class.
In \cite{mcc90c} it was shown that the locally convex inductive limit
topology on (\ref{eq:un}) was the Mackey topology 
of $(N^+ (B_1) ,\rho)$, and so a function $f$ is in the common
range of all co-analytic Toeplitz operators in $H^2 (B_1)$ if 
and only if  it is in the dual of $(N^+ (B_1) ,\rho)$;
the dual has been characterized by Yanagihara as those 
functions for which $\hat f(n) = O(e^{-c\sqrt n})$ for some $c > 0$
\cite{yan73a} (see \cite{mcc92a}
for a simpler proof).
\par
In \cite{naw89}, Nawrocki characterized the dual of $N^+(\bd)$ as
those functions $f$ whose Taylor coefficients at zero satisfy
\setcounter{equation}{\value{theorem}}
\begin{equation}
\hat f(\alpha) = O(e^{-c|\alpha|^{d/(d+1)}}).
\label{eq:na}
\end{equation}
\addtocounter{theorem}{1}
We prove that (\ref{eq:na}) does not characterize the common
range of co-analytic Toeplitz operators in $\hbd, d \geq 2$:
\par
{\bf Theorem \ref{thm:big}} {\it
Let $f(z_1,\dots,z_d) = f_1 (z_1) = \sum_{n=0}^\i a_n z_1^n$,
let $\varepsilon > 0$,
and suppose that $a_n = O(e^{-c n^{\frac 12 + \varepsilon}})$ for
some $c>0$. Then $f$ is in the range of the  Toeplitz
operator $\tmh$ for every non-zero $m$ in $\hibd$.
}
\vskip 5pt
It follows, therefore, that for $d \geq 2$, the functional
induced by 
$f(z_1,\dots,z_d) = \sum_{k=1}^\i e^{-k^{\frac 47 }} z_1^k $
is not bounded on some $P^2(w \sigma_d)$. 
In the proof of Theorem~\ref{thm:main} we construct
such a $w$ that is continuous and vanishes at only one point. 

\section{Preliminary Lemmata}
We need to know explicitly the projection from $L^2(\bd)$ onto
$H^2(\bd)$. Let  $\alpha= (\alpha_1,
\dots,\alpha_d)$ be a multi-index and $\zeta = (z_1,\dots,z_d)$
a point in $\cd$. The function $\zeta^\alpha$ then maps
$\zeta$ to $z_1^{\alpha_1} \dots z_d^{\alpha_d}$.
The notation $|\alpha|$ stands for $\alpha_1 + \dots \alpha_d$,
and $\alpha! = \alpha_1!\, \dots \alpha_d!$.
\begin{lemma}
\label{lem:co}
\setcounter{equation}{\value{theorem}}
\begin{equation}
\int_{S_d }  \zeta^\alpha \overline{\zeta^\beta} d \sigma_d
= \delta_{\alpha,\beta} \frac{ (d-1)!\, \alpha !}{(d-1 +|\alpha|)!} .
\label{eq:no}
\end{equation}
\addtocounter{theorem}{1}
Moreover, if $P_{H^2({B_d})}$ denotes the projection from
$L^2(\sigma_d)$ onto $H^2({B_d})$, then
\setcounter{equation}{\value{theorem}} 
\begin{equation} 
P_{H^2({B_d})} |z_2^{\alpha_2}|^2 \dots
|z_d^{\alpha_d}|^2 \overline{z_1^j} z_1^i =
\cases{0 &if $i < j$\cr
       \frac{(d-1+i-j)!\, i!\, \alpha_2 !\, \dots \alpha_d !}
             {(i-j)!\, (d-1 + i + \alpha_2 + \dots + \alpha_d)!}
               z_1^{i-j} 
          &if $i \geq j$.\cr}
\label{eq:pr} 
\addtocounter{theorem}{1}
\end{equation}
\end{lemma}
\vskip 5pt
{\sc Proof:} Formula (\ref{eq:no}) is proved in \cite{rud80}.
The expression on the left-hand side of (\ref{eq:pr})
is orthogonal to every monomial except $z_1^{i-j}$; taking
inner products gives the constant.
\hfill $\Box$
\par
\vskip 5pt
We wish to be able to transfer information about co-analytic
Toeplitz operators with the same symbol on different spaces.
To do this we use the following lemma, whose proof is immediate:
\begin{lemma}
\label{lem:to}
Let $\H$ be a Hilbert space of holomorphic functions on 
${B_d}$ in which the monomials are mutually orthogonal.
Let $m(z_1,\dots,z_d) = \sum_{\beta \in {\N}^d} b_{\beta} \zeta^\beta$.
Then 
\setcounter{equation}{\value{theorem}}
\begin{equation}
T^{\H}_{\bar m} \frac {\zeta^\alpha}{\V \zeta^\alpha \V^2_{\H}}
= \sum_{\beta \leq \alpha} \bar b_{\alpha - \beta}
\frac {\zeta^\beta}{\V \zeta^\beta \V^2_{\H}}
\label{eq:to} 
\addtocounter{theorem}{1}
\end{equation}
\end{lemma}
\vskip 5pt
This lemma also allows us to define Toeplitz operators
with an unbounded conjugate analytic symbol. The formal 
definition(\ref{eq:to}) defines an upper triangular
operator, with respect to the orthonormal basis of normalized monomials.
It therefore has a domain which contains all the polynomials;
we extend its domain to include all functions on which $T_{\bar m}$,
thought of as a formal operator on the power series, converges
in each entry to the Taylor coefficients of some function in $\H$.
\par
Let $A^{-n}$ consist of all holomorphic functions $m$ in the unit
disk, satisfying $|m(z)| = O( (1-|z|)^{-n})$. 
The space $A^0$ is $H^\i(B_1)$.
For $n=0$, the following result is proved in \cite{mcc93b}.
\begin{lemma}
\label{lem:ko}
Let $f$ be in $A^{-n}$ for some $n$, and $ 0 < \alpha < 2$.
Then
$$
\int_{B_1} ( \log^- |f|)^\alpha dA < \i .
$$
\end{lemma}
\vskip 5pt
{\sc Proof:}
We can assume that $f$ has no zeroes in $\{ z : |z| < \frac 12 \}$.
It follows from \cite{mcc93b} and standard Nevanlinna theory 
that for any function
$g$ in the Nevanlinna class $N(B_1)$, with $g(0) \neq 0$,
and any $a < 2$,
$$
\int_{B_1} ( \log^- |g|)^a dA 
\leq K( T(g,1), |g(0)|, a) .
$$
Fix $p$ strictly between $1$ and $\frac {2}{\alpha}$,
let $a = \alpha p < 2$,
let $q = \frac{p}{p-1}$ and let $N > q$.
Let $D_{1}$ be a smoothly bounded convex domain inside the disk,
containing $\{ z : |z| < \frac 12 \}$, 
whose closure touches the unit circle only at $1$,
and which has a high degree of tangency at $1$:
let the boundary of $D_1$ be $\{ \rho(\theta) e^{i\theta} : -\pi \leq
\theta \leq \pi \}$,
and assume $1 - \rho(\theta) \sim |\theta |^N$.
For any other point $\zeta = e^{i\theta_0}$ on
the boundary of the unit disk,
let $D_\zeta = e^{i \theta_0} D_1$.
\par
Let $\psi_\zeta$ be the Riemann map of $D_\zeta$
onto $B_1$ that takes $0$ to $0$ and $\zeta$ to $\zeta$.
As the boundary of $D_\zeta$ is smooth, it follows from
the Kellog-Warschawski theorem (see e.g. \cite{Pom92}) that $\psi_\zeta$ and
its derivatives extend continuously to the closure of
$D_\zeta$, so distances before and after the conformal mapping
are comparable.
\par
If $r < \frac 1N$, then $f$ is in $H^r(D_\zeta)$,
and $\sup_{\zeta \in S_1} \V f \circ \psi_\zeta^{-1} \V_{H^r} < \i$.
Thus
$$
\int_{D_{e^{i\theta}}} (\log^- |f|)^a dA \leq C, \quad for\ all\ 
e^{i\theta} .
$$
Integrating with respect to $\theta$ and changing the order of
integration yields
$$
\int_{B_1} (\log^-|f(re^{i\phi})|)^a (1-r)^{\frac 1N} r dr d\phi
< \i .
$$
Now
\begin{eqnarray*}
\int_{B_1} ( \log^- |f|)^\alpha dA  \leq &
\left[ \int_{B_1} (\log^- |f|)^{\alpha p} (1-r)^{\frac pN} dA \right]^{\frac 1p}
\left[ \int_{B_1} (1-r)^{- \frac qN} dA \right]^{\frac 1 q} \\
< & \i .
\end{eqnarray*}
\hfill $\Box$
\par
\vskip 5pt
Let $\mu_n$ be
the measure on the unit disk given by
$d\mu_n(z) = \frac{1}{\pi} (1-|z|^2)^n dArea(z)$, and let ${\H}_n$ be 
$P^2(\mu_n)$. It is routine to verify that in ${\H}_n$  the
monomials are mutually orthogonal, and 
$$
\V z^k \V_{{\H}_n}^2 = \frac{1}{(k+1) \dots (k+n+1)}.
$$
The space ${\H}_0$ is the usual Bergman space for the disk.
The following lemma is proved in \cite{mcc93b} for
$n =0$;
for $ n > 0$, basically the same proof works (given
Lemma~\ref{lem:ko}), though
some care must be taken as $\tm$ is no longer bounded.
We include a proof for completeness.
\begin{lemma}
\label{lem:hn}
Let $n \geq 0$, and $m$ be a function in  $A^{-n}$, not
identically zero. Suppose $f(z) = \sum_{k=0}^\i a_k z^k$
where $a_k = O(e^{-c k^{\frac 12 + \varepsilon}})$ for
some $\varepsilon$ and $c$ greater than $0$. 
Then for any $s\geq 2n$ there exists $g$ in $\hs$ such that
$T^{\hs}_{\bar m}  g = f$.
\end{lemma}
\vskip 5pt
{\sc Proof:}
First, observe that $f = T^{\hs}_{\bar m}  g$  for
some $g$ if and only if
there is a constant $C$ such that for all polynomials $p$
$$
| \la p, f \ra_{\hs} | \leq C \sqrt{ \int |p|^2 |m|^2 d\mu_{s} }.
$$
So it is sufficient to prove that 
$$
| \sum_{k=0}^\i \bar a_k \hat p (k) 
\frac {1}{(k+1) \dots (k + s +1)} | 
\leq C \sqrt{ \int |p|^2 |m|^2 d\mu_{s}} .
$$
This in turn will follow from the Banach-Steinhaus
theorem if we can show that for any function $h$ in
$P^2(|m|^2 \mu_{s})$, 
\setcounter{equation}{\value{theorem}}
\begin{equation}
\hat h (k) = O(e^{ c k^{\frac 12 + \varepsilon}}).
\label{eq:grz}
\end{equation}
\addtocounter{theorem}{1}
Now Stoll showed in \cite{sto77} that 
if $h$ satisfies
$$
\int_{B_1} (\log^+ |h|)^\alpha d Area < \i 
$$
for some $\alpha > 0$ then $\hat h (k) = O(e^{o(\frac{2}{2+\alpha})})$.
We can assume $\varepsilon$ is small, and take
$\alpha = \frac{2-4\varepsilon}{1+2\varepsilon}$.
As $h$ is in $P^2(|m|^2 \mu_{s})$,
$h(z) m(z) (1-|z|^2)^{s/2} k(z)$ is in $L^2(dArea)$,
and 
$$
\log^+|h| \leq \log^+|k| + \log^-|(1-|z|^2)^{s/2}|
+ \log^-|m| .
$$
The first
two terms on the right are clearly integrable to the $\alpha^{th}$
power, and so is 
the third by Lemma~\ref{lem:ko};
therefore $h$ satisfies (\ref{eq:grz}) as desired.
\hfill $\Box$
\vskip 5pt
\par
We want to be able to restrict functions in the ball to 
planes and factor out zeros without losing control
of the size of the function; the next lemma allows us to do this.
\begin{lemma}
\label{lem:dim}
Let $m$ be holomorphic on $B_d$ and satisfy
$$
|m(z_1, \dots,z_d )| \leq C (1 - \sqrt{|z_1|^2 + \dots
|z_d|^2} )^{-s} .
$$
Suppose also that 
$$
m(z_1, \dots , z_d ) = z_d^t m_2(z_1, \dots , z_d) + z_d^{t+1}
m_3(z_1, \dots , z_d) ,
$$
where $m_2$ and $m_3$ are analytic. Let $m_1 (z_1, \dots, z_{d-1})
= m_2(z_1, \dots , z_{d-1}, 0)$. 
Then 
$$
|m_1 (z_1, \dots, z_{d-1})| \leq (3d)^{s+t} C (1 - \sqrt{|z_1|^2 + \dots
|z_{d-1}|^2})^{-(s+t)} .
$$
\end{lemma}
\vskip 5pt
{\sc Proof:}
Let $(z_1,\dots,z_{d-1})$ be in $B_{d-1}$, and let $\varepsilon =
\frac{1}{3d} ( 1 - \sqrt{|z_1|^2 + \dots + |z_{d-1}|^2}$).
Then the polydisk centered at $(z_1,\dots,z_{d-1},0)$ with multi-radius
$(\varepsilon,\dots,\varepsilon)$ is contained in $(
1-\varepsilon)B_d$.
Integrating on the distinguished boundary of the polydisk
we get
\begin{eqnarray*}
| m_1 (z_1,\dots,z_{d-1}) | = & | m_2 (z_1,\dots,z_{d-1},0) | \\ 
=& \left\vert {\int \atop (z_1,\dots,z_{d-1},0) + \varepsilon T^d}
\frac{m(\zeta_1, \dots,\zeta_d)}{\zeta_3^t} \right\vert
\\
\leq& \frac{C}{\varepsilon^{s+t}}.
\end{eqnarray*}
\hfill $\Box$
\section{Common Range of $T_{\bar m}$}
We can now prove that a function that depends on only one
variable is in the range of every $T^{\hbd}_{\bar m}$ if its
Taylor coefficients decay like $e^{-c k^{\frac{1}{2} + \varepsilon}}$.
\begin{theorem}
\label{thm:big}
Let $f(z_1,\dots,z_d) = f_1 (z_1) = \sum_{n=0}^\i a_n z_1^n$,
and suppose that $a_n = O(e^{-c n^{\frac 12 + \varepsilon}})$ for
some $\varepsilon,c>0$. Then $f$ is in the range of the  Toeplitz
operator $\tmh$ for every non-zero $m$ in $\hibd$.
\end{theorem}
\vskip5pt
{\sc Proof:} For $d =1$, this is proved (without the
$\varepsilon$) in \cite{mcc90c}, so assume $d \geq 2$. Fix $m$ in $\hibd$;
$$
m \zs = \sum_{i_1,\dots,i_d =0}^\i b_{i_1,\dots,i_d} 
z_1^{i_1} \dots z_d^{i_d}  .
$$
Let 
$$
S = \{ (i_2,\dots, i_d ): {\rm for\ some\ }i_1,\, b_{i_1,\dots,i_d} \neq 0 \}.
$$
Define
$$
t_d = \inf\{i_d: {\rm for\ some\ }i_2,\dots,i_{d-1},\, 
(i_2,\dots,i_{d-1},i_d) \in S \},
$$
and define $t_k$ inductively by
$$
t_k = \inf\{i_k: {\rm for\ some\ }i_2,\dots,i_{k-1},\, 
(i_2,\dots,i_{k-1},i_k,t_{k+1},\dots,t_d) \in S \}.
$$
Let $n = t_2 + \dots + t_d$.
\par
Case (a): $n = 0$.
\par
Then the function
$$
m_1(z_1) = m(z_1,0,\dots,0)
$$
is not identically zero, and is in $\hid$.
By Lemma~\ref{lem:co}, 
$$
\tmh z_1^i =
\sum_{j} \bar b_{j,0,\dots, 0} \frac{(i-j+1) \dots (i-j+d-1)}{(i+1)
\dots (i+d-1)} z_1^{i-l} .
$$
So by Lemma~\ref{lem:to}, if one can solve
the equation
\setcounter{equation}{\value{theorem}}
\begin{equation}
T_{\bar m_1}^{{\H}_{d-2}} g_1 = f_1
\label{eq:r1}
\end{equation}
\addtocounter{theorem}{1}
for some $g_1$ in ${\H}_{d-2}$, then $g \zs = g_1 (z_1)$ solves
$$
\tmh g = f ,
$$
and by Equation~(\ref{eq:no})
$\V g \V_{\hbd} = \sqrt{(d-1)!} \V g_1 \V_{{\H}_{d-1}} < \i$.
By Lemma~\ref{lem:hn}, equation~(\ref{eq:r1}) has a 
solution.
\par
Case (b): $n > 0$.
\par
One can decompose $m$ as
$$
m \zs = z_2^{t_2} \dots z_d^{t_d} m_2 \zs +  
 m_3 \zs ,
$$
where each term in the expansion of $m_3$ is divisible
by some $z_k^{t_k + 1}$. 
Applying Lemma~\ref{lem:dim} inductively,
$m_1(z) = m_2(z,0,\dots,0)$ is in $A^{-n}$,
and by the choice of $t_2, \dots, t_d$, it is not identically zero.
Consider the function 
$$
f_2(z) = \sum_{k=0}^\i a_k (k+d) (k+d+1) \dots (k+ dn+1) z^k.
$$
As $d \geq 2$, we can apply
Lemma~\ref{lem:hn} with $s = dn$, 
so there is 
$$
g_2(z) = \sum_{k=0}^\i \gamma_k (k+1) (k+2) \dots (k+ dn+1) z^k
$$
in ${\cal H}_{dn}$ with 
\setcounter{equation}{\value{theorem}}
\begin{equation}
T_{\bar m_1}^{{\cal H}_{dn}} g_2 = f_2. 
\label{eq:r2}
\end{equation}
\addtocounter{theorem}{1}
Define $g$ by
$$
g\zs = \frac{1}{t_2!\, \dots t_d!} z_2^{t_2}\dots z_d^{t_d} \sum_{k=0}^\i 
\gamma_k (k+1) (k+2) \dots (k+ n+d-1) z_1^k.
$$
The function $g$ is in $\hbd$ because
\begin{eqnarray*}
\V g \V^2_{\hbd} =& \frac{(d-1)!}{t_2!\, \dots t_d!} 
\sum_{k=0}^\i | \gamma_k |^2 (k+1) \dots (k+n+d-1) \\
\leq &\frac{(d-1)!}{t_2!\, \dots t_d!}
\sum_{k=0}^\i | \gamma_k |^2 (k+1) \dots (k+dn+1)\\
=&  \frac{(d-1)!}{t_2!\, \dots t_d!}
\V g_2 \V^2_{{\cal H}_{(d-1)n}} \\
<& \i.
\end{eqnarray*}
Moreover 
$$
\tmh g = T_{\bar z_2^{t_d} \dots \bar z_d^{t_d} 
\overline{m_1(z_1)}}^{H^2({\bd})} g
$$
is a function of $z_1$ only; it is, in fact, $f$.
For if $\tmh g = \sum_{k=0}^\i e_k z_1^k$, and
$m_1(z) = \sum_{k=0}^\i c_k z^k$, then 
taking the inner product with $z_1^j$ we get 
\begin{eqnarray}
\frac{(d-1)!}{(j+1)\dots (j+d-1)} e_j =&
 \la \tmh g , z_1^j \ra_{\hbd} \nonumber \\
=& \la g, z_2^{t_2} \dots  z_d^{t_d} m_1  z_1^j \ra_{\hbd} \nonumber \\
=&(d-1)! \sum_{k=j}^\i \gamma_k  \bar c_{k-j} 
\label{eq:r3}
\addtocounter{theorem}{1}
\end{eqnarray}
Taking the inner product with $z^j$ in Equation(\ref{eq:r2}),
we get
\begin{eqnarray}
\frac{1}{(j+1)\dots (j+d-1)} a_j =&
\la T^{{\cal H}_{dn}}_{\bar m_1} g_2 , 
z^j \ra_{{\cal H}_{dn}} \nonumber \\
=& \la g_2, m_1 z^j \ra_{{\cal H}_{dn}} \nonumber \\
=& \sum_{k=j}^\i \gamma_k  \bar c_{k-j}
\label{eq:r4}
\addtocounter{theorem}{1}
\end{eqnarray}
Comparing Equations(\ref{eq:r3}) and (\ref{eq:r4}), we see that
$\tmh g = f$, as desired.
\hfill $\Box$
\vskip 5pt
\par
\section{Proof of the main theorem}
Define $\fc$ by
\setcounter{equation}{\value{theorem}}
\begin{equation}
\fc(z) = \exp(c \frac {1 - |w|^2}{(1- \la z, w \ra)^{d+1}} )
\label{eq:f}
\end{equation}
\addtocounter{theorem}{1}
We need the following two results. The first
states that $d(c\fc , 0) \to 0$, and was proved by Drewnowski.
A proof is given in \cite[Lemma 3.2]{naw89}.
\begin{lemma}[Drewnowski]
\label{lem:log}
$$
\lim_{c \to 0} \quad \sup_{w \in \bd} \int_{S_d}
 \log( 1 + |c \fc |) d \sigma_d = 0 .
$$
\end{lemma}
The second result, due to Nawrocki, estimates the growth of the Taylor
coefficients of $\fc$. We are interested in $w = re_1 = 
(r,0,\dots,0)$; in this case all the Taylor coefficients
of $F_{c,{re_1}}$ are positive, and the following follows
easily from the proof of \cite[Lemma 3.3]{naw89}:
\begin{lemma}[Nawrocki]
\label{lem:tay}
For each $c > 0$ there exists $\varepsilon > 0$
such that
$$
\inf_{i \in \N} \sup_{ 0 < r < 1}
\sqrt{\frac{(d-1+i)!}{(d-1)!\, i!} }
\hat F_{c,re_1} (i,0,\dots,0)
e^{-\varepsilon i^{\frac{d}{d+1}}}
> 0 .
$$
\end{lemma}
We can now prove the main theorem.
\begin{theorem}
\label{thm:main}
Let $d \geq 2$. There is a continuous non-negative function $g$
on $S_d$, vanishing only at the point $e_1$, 
and  satisfying 
$\int_{S_d} log(g) d\sigma_d > -\i$, with the property that
the only function $m$ in $\hibd$ with $|m^\star| \leq g$ almost
everywhere $[\sigma_d]$ is the zero function.
\end{theorem}
\vskip5pt
{\sc Proof:}
Let $V_n = \{ \zeta \in S_d : | \zeta - e_1 | \geq \frac 1n \}$.
By Lemma~\ref{lem:tay},  for any sequence $c_n$ tending to zero,
one can choose $i_n$ and $r_n$ such that
\setcounter{equation}{\value{theorem}}
\begin{equation}
\hat F_{c_n,r_n e_1} (i_n, 0, \dots, 0) > 
\frac{n}{c_n}
e^{(i_n)^{\frac 47}} 
\label{eq:jz}
\end{equation}
\addtocounter{theorem}{1}
(because $\frac 4 7 < \frac {d}{d+1}$).
Moreover, by passing to a subsequence, one can assume that
\setcounter{equation}{\value{theorem}}
\begin{equation}
\sup_{\zeta \in V_n} c_n |F_{c_n,r_n e_1}(\zeta)| \leq \frac{1}{2^n} ,
\label{eq:su}
\end{equation}
\addtocounter{theorem}{1}
and that
$$
\int_{S_d} \log( 1 + |c_n F_{c_n,r_n e_1} |) d \sigma_d \leq \frac{1}{2^n}.
$$
Define $g$ by
$$
g(\zeta) = \sqrt{ \frac{1}{1+ \sum_{n=1}^\i |c_n F_{c_n, r_n e_1} (\zeta) |^2}} .
$$
It follows from (\ref{eq:su}) that $g$ is continuous and vanishes only at
$e_1$. Moreover
\begin{eqnarray*}
\int_{S_d} \log g d\sigma_d =& - \frac 12 
\int_{S_d} \log(1+ \sum_{n=1}^\i |c_n F_{c_n, r_n e_1}|^2) d\sigma_d \\
> & - \int_{S_d} \log \prod_{n=1}^\i (1 + |c_n F_{c_n, r_n e_1}|)^2 d\sigma_d \\
= & - 2 \sum_{n=1}^\i \int_{S_d} \log( 1 + |c_n F_{c_n,r_n e_1} |) d \sigma_d \\
\geq & -2 .
\end{eqnarray*}
Now suppose there is a non-zero $m$ in $\hibd$ with $|m^\star | \leq g$ a.e.
Then each of the functions $c_n F_{c_n, r_n e_1}$, being analytic
in the ball of radius $\frac{1}{r_n}$, is in $P^2(|m^\star |^2 \sigma)$;
moreover they are all of norm less than one in this space,
because 
$$
\int_{S_d} | c_n F_{c_n, r_n e_1} |^2 |m^\star |^2 d\sigma
\leq \int_{S_d} | c_n F_{c_n, r_n e_1} |^2  g^2 d\sigma < 1 .
$$
\par
Let 
$$
f \zs = \sum_{k=0}^\i e^{-k^{4 \over 7}} \frac{(k+d-1)!}{(d-1)!k!} z_1^k .
$$
By Theorem~\ref{thm:big}, there is a function $h$ in $\hbd$
with 
$$
\tmh h = f .
$$
It follows that the linear map
$$
\Gamma: p \mapsto \la p , f \ra_{\hbd} ,
$$
defined a priori on the polynomials, extends by continuity
to a bounded linear map on $P^2(|m^\star |^2 \sigma)$,
as
$$
| \Gamma(p) | =|\la p,P(\bar m h) \ra| =|\int p m^\star \bar h^\star
d\sigma_d | \leq \V h \V_{\hbd} \V p \V_{P^2(|m^\star |^2 \sigma)}
.
$$
Moreover, 
each function $c_n F_{c_n, r_n e_1}$ is
uniformly approximated on $S_d$ by
the partial sums of its Taylor series;
hence
\setcounter{equation}{\value{theorem}}
\begin{equation}
\Gamma(c_n F_{c_n, r_n e_1}) = \sum_{k=0}^\i
c_n \hat F_{c_n, r_n e_1} (k) e^{-k^{4 \over 7}}.
\label{eq:iz}
\end{equation}
\addtocounter{theorem}{1}
But all the terms on the right-hand side of (\ref{eq:iz})
are positive, and the $i_n^{th}$ term is at least $n$ by Equation~(\ref{eq:jz}).
This contradicts the boundedness of $\Gamma$.
\hfill $\Box$
\par
We note that the theorem is much easier to prove if $g$
is not required to be continuous.
One shows, as in the one-variable case, that
$$
\nb = \cup P^2(w \sigma_d) 
$$
where $w$ ranges over all log-integrable weights.
This is a strictly larger set than
$$
\cup_{m \in \hibd} P^2(|m^\star |^2 \sigma_d)
$$
because if $f$ is in some $P^2(|m^\star |^2 \sigma_d)$,
then $fm$ is in $\hbd$, so the zero-set of
$f$ is contained in an $\hbd$-zero set; 
but it is a result of Rudin
that for any $p < 2$, there is a function in
$H^p(\bd)$ whose zero set is not contained in any
$\hbd$-zero set \cite{rud76}. It follows that such a function
is in $P^2(w^2 \sigma_d)$ for some log-integrable $w$,
but not in any $P^2 (|m^\star |^2 \sigma_d)$, so $w$
cannot dominate the modulus of any $m^\star$.


\begin{thebibliography}{10}

\bibitem{mcc93b}
J.E. M\raise.5ex\hbox{c}Carthy.
\newblock Coefficient estimates in weighted {B}ergman spaces.
\newblock to appear.

\bibitem{mcc90c}
J.E. M\raise.5ex\hbox{c}Carthy.
\newblock Common range of co-analytic {T}oeplitz operators.
\newblock {\em J. Amer. Math. Soc.}, 3(4):793--799, 1990.

\bibitem{mcc92a}
J.E. M\raise.5ex\hbox{c}Carthy.
\newblock Topologies on the {Smirnov} class.
\newblock {\em J. Funct. Anal.}, 104(1):229--241, 1992.

\bibitem{naw89}
M.~Nawrocki.
\newblock Linear functionals on the {S}mirnov class of the unit ball in
  {$C^n$}.
\newblock {\em Annales Acad. Sci. Fenn.}, 14:369--379, 1989.

\bibitem{Pom92}
C.~Pommerenke.
\newblock {\em Boundary behaviour of conformal maps}.
\newblock Springer-Verlag, Berlin, 1992.

\bibitem{rud76}
W.~Rudin.
\newblock Zeros of holomorphic functions in balls.
\newblock {\em Indag. Math.}, 38:57--65, 1976.

\bibitem{rud80}
W.~Rudin.
\newblock {\em Function Theory in the unit ball of $C^n$}.
\newblock Springer-Verlag, Berlin, 1980.

\bibitem{rud86}
W.~Rudin.
\newblock {\em New constructions of functions holomorphic in the unit ball of
  $C^n$}.
\newblock C.B.M.S.\ No.~63. American Mathematical Society, Providence, 1986.

\bibitem{sto77}
M.~Stoll.
\newblock Mean growth and {T}aylor coefficients of some topological algebras of
  analytic functions.
\newblock {\em Ann. Polon. Math.}, XXXV:139--158, 1977.

\bibitem{sze21}
G.~{Szeg\"o}.
\newblock {\"Uber} die {R}andwerten einer analytischen {F}unktionen.
\newblock {\em Math. Ann.}, 84:232--244, 1921.

\bibitem{yan73a}
N.~Yanagihara.
\newblock Multipliers and linear functionals for the class {$N^+$}.
\newblock {\em Trans. Amer. Math. Soc.}, 180:449--461, 1973.

\end{thebibliography}


\end{document}